\documentclass[12pt, a4paper]{article}
\usepackage{amssymb}
\usepackage{latexsym}
\usepackage{amsbsy}

\newcommand{\bcen}{\begin{center}}     \newcommand{\ecen}{\end{center}}
\newcommand{\bay}{\begin{array}}      \newcommand{\eay}{\end{array}}
\newcommand{\beq}{\begin{eqnarray*}}      \newcommand{\eeq}{\end{eqnarray*}}
\def\az{\alpha}
\def\bz{\beta}

\def\ann{\mbox{ann}}
\def\Aut{\mbox{Aut}}

\def\char{\mbox{char}}

\def\gr{\mbox{gr}}
\def\Gr{\mbox{Gr}}
\def\Hom{\mbox{Hom}}
\def\HOM{\mbox{HOM}}

\def\Ext{\mbox{Ext}}
\def\EXT{\mbox{EXT}}
\def\End{\mbox{End}}

\def\dim{\mbox{dim}}
\def\mod{\mbox{mod}}
\def\Mod{\mbox{Mod}}

\def\Im{\mbox{Im}}

\begin{document}

\title{Construction of Koszul algebras by finite Galois covering \footnote{Project 10201004 supported by NSFC.}}

\author{Yang Han $^a$ and Deke Zhao $^b$}

\date{\footnotesize a. Academy of Mathematics and Systems Science,
Chinese Academy of Sciences,\\ Beijing 100080, P.R.China. E-mail:
hany@iss.ac.cn\\ b. School of Mathematical Science, Beijing Normal
University, Beijing 100875,\\ P.R. China. E-mail:
dkzhao@mail.bnu.edu.cn}

\maketitle

\begin{abstract}

It is shown that, the quasi-Koszulities of algebras and modules are
Morita invariance. A finite-dimensional $K$-algebra $A$ with an
action of $G$ is quasi-Koszul if and only if so is the skew group
algebra $A \ast G$, where $G$ is a finite group satisfying $\char K
\nmid |G|$. A finite-dimensional $G$-graded $K$-algebra $A$ is
quasi-Koszul if and only if so is the smash product $A \# G^*$,
where $G$ is a finite group satisfying $\char K \nmid |G|$. These
results are applied to prove that, if a finite-dimensional connected
quiver algebra is Koszul then so are its Galois covering algebras
with finite Galois group $G$ satisfying $\char K \nmid |G|$. So one
can construct Koszul algebras by finite Galois covering. Moreover, a
general construction of Koszul algebras by Galois covering with
finite cyclic Galois group is provided. As examples, many Koszul
algebras are constructed from exterior algebras and Koszul
preprojective algebras by finite Galois covering with either cyclic
or noncyclic Galois group.

\end{abstract}


\section*{Introduction}

Koszul algebras play an important role in commutative algebra,
algebraic geometry, algebraic topology, Lie theory and quantum
groups (cf. \cite{P, BF, L, BGSc, Bondal, Rosso, BGS}). It is a
quite nice class of algebras because, on one hand, there exists
Koszul duality in the sense of not only algebra (cf. \cite[Theorem
2.10.2]{BGS} and \cite[Theorem 2.3]{GM1998}) but also module
category and derived category (cf. \cite[Theorem 5.2]{GM1998} and
\cite[Theorem 2.12.6]{BGS}), on the other hand, for a Koszul
algebra, both the minimal projective resolution of its semisimple
part and its minimal projective bimodule resolution are rather clear
(cf. \cite{BK} and \cite{GHMS}). Now many algebras are known to be
Koszul, such as path algebras \cite[p. 240]{GM1998}, radical square
zero algebras \cite[Corollary]{Mp}, quadratic global dimension 2
algebras \cite[Theorem 7.2]{GM1996}, finite-dimensional
indecomposable radical cube zero selfinjective algebras of infinite
representation type \cite[Theorem 1.5]{Mp} and many (generalized)
preprojective algebras (cf. \cite[Section 7]{G1997} and
\cite[Theorem 1.9]{Mp}). From the known Koszul algebras, there are
some ways to construct new Koszul algebras: the opposite algebra of
a Koszul algebra is Koszul (cf. \cite[Proposition 2.2.1]{BGS} and
\cite[Corollary 4.3]{GM1998}), the quadratic duality equivalently
the Yoneda algebra of a Koszul algebra is Koszul (cf.
\cite[Proposition 2.9.1]{BGS}, \cite[Theorem 6.1]{GM1996} and
\cite[Theorem 2.2]{GM1998}), the tensor product algebra of two
Koszul algebras is Koszul \cite[Theorem 3.7]{GM1998}, and so on.

Here we provide another construction of Koszul algebras by finite
Galois covering. This construction is very convenient in practice
because one can construct Koszul algebras and provide their quivers
and relations directly from a known Koszul algebra given by quiver
with relations. The paper is organized as follows: In
Section~\ref{moritaequivalence}, we show that the quasi-Koszulities
of algebras and modules are Morita invariance.
Section~\ref{skewgroupalgebra} is essentially due to
Mart\'{i}nez-Villa (cf. \cite{Ms}). In this section, we show that a
finite-dimensional $K$-algebra $A$ with an action of $G$ is
quasi-Koszul if and only if so is the skew group algebra $A \ast G$,
where $G$ is a finite group satisfying $\char K \nmid |G|$.
Moreover, we prove that the Yoneda algebra $\Ext^*_{A \ast G}((A
\ast G)/J(A \ast G), (A \ast G)/J(A \ast G))$ of the skew group
algebra $A \ast G$ is isomorphic to the skew group algebra
$\Ext^*_A(A/J,A/J) \ast G$ of the Yoneda algebra
$\Ext^*_A(A/J,A/J)$. In Section~\ref{smashproduct}, we show in two
ways that a finite-dimensional $G$-graded $K$-algebra $A$ is
quasi-Koszul if and only if so is the smash product $A \# G^*$,
where $G$ is a finite group satisfying $\char K \nmid |G|$.
Moreover, we prove that the Yoneda algebra $\Ext^*_{A \# G^*}((A \#
G^*)/J(A \# G^*), (A \# G^*)/J(A \# G^*))$ of the smash product $A
\# G^*$ is isomorphic to the smash product $\Ext^*_A(A/J, A/J) \#
G^*$ of the Yoneda algebra $\Ext^*_A(A/J,A/J)$. In
Section~\ref{galoiscovering}, we apply the results above to show
that, if a finite-dimensional connected quiver algebra is Koszul
then so are its Galois covering algebras with finite Galois group
$G$ satisfying $\char K \nmid |G|$. Finally, in
Section~\ref{construction}, we provide a general construction of
Koszul algebras by Galois covering with finite cyclic Galois group.
As examples, we construct many Koszul algebras from exterior
algebras and Koszul preprojective algebras by finite Galois covering
with either cyclic or noncyclic Galois group.

\section{Morita invariance of quasi-Koszulity} \label{moritaequivalence}

\newtheorem{theorem}{Theorem}
\newtheorem{proposition}{Proposition}
\newtheorem{lemma}{Lemma}
\newtheorem{corollary}{Corollary}
\newtheorem{remark}{Remark}
\newtheorem{example}{Example}

In this section we show that the quasi-Koszulities of algebras and
modules are Morita invariance. Throughout the paper, the composition
of maps is written from right to left except for those in Yoneda
algebras.

\subsection{Quasi-Koszulity}

Let $K$ be a field. A graded $K$-algebra $A=\coprod_{i \geq 0}A_i$
is said to be {\it generated in degree 0 and 1} if $A_i=A^i_1$ for
all $i \geq 2$ (cf. \cite{GM1998}).

In this section, $A$ is assumed to be a fixed noetherian semiperfect
$K$-algebra (cf. \cite[\S 27]{AF}). Denote by $J:=J(A)$ the Jacobson
radical of $A$. Let $M, N$ be two $A$-modules. Then $\Ext^i_A(M,N)$,
$i \geq 1$, can be viewed as a $K$-vector space of congruence
classes of $i$-extensions whose addition is given by Baer sum (cf.
\cite[Chapter III, \S 5]{Mac}). It is convenient to write an
$i$-extension in the form $M \twoheadleftarrow E_1 \leftarrow \cdots
\leftarrow E_i \leftarrowtail N$. It is well-known that
$\Ext^*_A(M,M) := \coprod_{i \geq 0}\Ext^i_A(M,M)$ is a graded
$K$-algebra with $\Ext^i_A(M,M)$ in degree $i$, for which the
multiplication is given by the Yoneda product (cf. \cite[Chapter
III, \S 5]{Mac}). The algebra $A$ is called a {\it quasi-Koszul
algebra} if $\Ext^*_A(A/J,A/J)$ is generated in degree 0 and 1 (cf.
\cite[p. 263]{GM1996}). An $A$-module $M$ is said to be {\it
quasi-Koszul} if $\Ext^i_{A}(M, A/J) \Ext^1_{A}(A/J,
A/J)=\Ext^{i+1}_{A}(M, A/J)$ for all $i \geq 0$, equivalently,
$\Ext^i_{A}(M, A/J) \Ext^j_{A}(A/J, A/J)=\Ext^{i+j}_{A}(M, A/J)$ for
all $i,j \geq 0$.

A graded $K$-algebra $A=\coprod_{i \geq 0}A_i$ is called a {\it
graded quiver algebra} if it satisfies: (1): $A_0 \cong K^r$ as
$K$-algebras for some $r \geq 1$; (2): $\dim_KA_i < \infty$ for all
$i \geq 0$; (3): $A$ is generated in degree 0 and 1 (cf. \cite{Mg}).
Note that a graded quiver algebra $A$ is isomorphic to a graded
quotient of a path algebra with length grading, namely $A \cong
KQ/I$ where $Q$ is a finite quiver and $I \subseteq (KQ^+)^2$ is a
homogeneous ideal of $A$ in the length grading, here $KQ^+$ denotes
the ideal of $KQ$ generated by all arrows (cf. \cite{G1997}). For
the theory of quivers and their representations we refer to
\cite{ARS}. A quasi-Koszul graded quiver algebra is called a {\it
Koszul algebra}.

\subsection{Decomposition of extension groups}\label{doeg}

Let $\oplus^s_{i=1}M_i$ be the direct sum of $A$-modules
$M_1,...,M_s$. Denote by $\lambda_i: M_i \rightarrow
\oplus^s_{i=1}M_i, m_i \mapsto (0, ... , 0, m_i, 0, ... , 0)$ and
$\rho_i : \oplus^s_{i=1}M_i \rightarrow M_i, (m_1,..., m_s)
\mapsto m_i$ the canonical injection and projection respectively.
Denote by $\Delta : M \rightarrow M^s, m \mapsto (m,...,m)$ and
$\nabla : M^s \rightarrow M, (m_1,...,m_s) \mapsto m_1+ \cdots
+m_s$ the {\it diagonal map} and the {\it sum map} respectively.
Note that $\Delta = \sum^s_{i=1} \lambda _i$ and $\nabla =
\sum^s_{i=1} \rho _i$.

For $r \geq 0$, $\Ext^r_A(\oplus^s_{i=1}M_i,N) \cong
\oplus^s_{i=1}\Ext^r_A(M_i,N)$ as $K$-vector spaces and the
isomorphism is given by $\phi : \Ext^r_A(\oplus^s_{i=1}M_i,N)
\rightarrow \oplus^s_{i=1} \Ext^r_A(M_i,N), [\xi] \mapsto
(\Ext^r_A(\lambda _i,N)([\xi]))_i$ with inverse $\psi :
\oplus^s_{i=1} \Ext^r_A(M_i,N) \rightarrow
\Ext^r_A(\oplus^s_{i=1}M_i,N),$ \linebreak $([\xi_i])_i \mapsto
\sum^s_{i=1}\Ext^r_A(\rho_i,N)([\xi_i]).$

Similarly, $\Ext^r_A(M,\oplus^t_{j=1}N_j) \cong
\oplus^t_{j=1}\Ext^r_A(M,N_j)$ as $K$-vector spaces and the
isomorphism is given by $\phi : \Ext^r_A(M,\oplus^t_{j=1}N_j)
\rightarrow \oplus^t_{j=1} \Ext^r_A(M,N_j), [\xi] \mapsto
(\Ext^r_A(M,\rho_j)([\xi]))_j$ with inverse $\psi :
\oplus^t_{j=1}\Ext^r_A(M,N_j) \rightarrow
\Ext^r_A(M,\oplus^t_{j=1}N_j),$ \linebreak $([\xi_j])_j \mapsto
\sum^t_{i=1}\Ext^r_A(M, \lambda_j)([\xi_j]).$

More general, $\Ext^r_A(\oplus^s_{i=1}M_i,\oplus^t_{j=1}N_j) \cong
\oplus^s_{i=1}\oplus^t_{j=1}\Ext^r_A(M_i,N_j)$ as $K$-vector
spaces and the isomorphism is given by $\phi :
\Ext^r_A(\oplus^s_{i=1}M_i,\oplus^t_{j=1}N_j) \rightarrow
\oplus^s_{i=1}\oplus^t_{j=1} \Ext^r_A(M_i,N_j), [\xi] \mapsto
(\Ext^r_A(M_i,\rho_j)\Ext^r_A(\lambda
_i,\oplus^t_{j=1}N_j)([\xi]))_{i,j}$ with inverse $\psi :
\oplus^s_{i=1}\oplus^t_{j=1} \Ext^r_A(M_i,N_j) \rightarrow
\Ext^r_A(\oplus^s_{i=1}M_i,\oplus^t_{j=1}N_j), ([\xi_{ij}])_{i,j}
\mapsto \sum^s_{i=1}\sum^t_{j=1}\Ext^r_A(\rho_i,
\oplus^t_{j=1}N_j) \Ext^r_A(M_i,\lambda_j)([\xi_{ij}]).$ In this
case, denote by \linebreak $E^r_A(M_i,N_j)$ the image of
$\Ext^r_A(M_i,N_j)$ under the map $\psi$. Then \linebreak
$\Ext^r_A(\oplus^s_{i=1}M_i,\oplus^t_{j=1}N_j) =
\coprod^s_{i=1}\coprod^t_{j=1}E^r_A(M_i,N_j)$.

\subsection{Morita invariance of quasi-Koszulity}

Firstly, we have the following ``expansion lemma'':

\begin{lemma} \label{sumlemma} Let $M, N, L_i$ be
$A$-modules with $1 \leq i \leq z$. Then, for $r , l \geq 0$, as the
$K$-subspaces of $\Ext^{r+l}_A(M,N)$,
$$\Ext^r_A(M, \oplus ^z_{i=1}L_i)\Ext^l_A(\oplus ^z_{i=1}L_i, N) =
\sum^z_{i=1}\Ext^r_A(M, L_i)\Ext^l_A(L_i, N).$$ \end{lemma}

{\bf Proof.} For $[\xi] \in \Ext^r_A(M, \oplus ^z_{i=1}L_i)$ and
$[\zeta] \in \Ext^l_A(\oplus ^z_{i=1}L_i, N)$, we have
$$\begin{array}{lll} &&
[\xi]\Ext^l_A(\lambda_i\rho_i,N)([\zeta])\\ & = & \Ext^r_A(M,
\rho_i)([\xi])\Ext^l_A(\lambda_i,N)([\zeta])\\ & = & \Ext^r_A(M,
\lambda_i\rho_i)([\xi]) [\zeta]. \end{array}$$ Thus
$$\begin{array}{lll} [\xi][\zeta] & = &
[\xi]\Ext^l_A(\sum^z_{i=1}\lambda_i\rho_i,N)([\zeta])\\ & = &
\sum^z_{i=1} [\xi]\Ext^l_A(\lambda_i\rho_i,N)([\zeta])\\ & = &
\sum^z_{i=1} \Ext^r_A(M,
\rho_i)([\xi])\Ext^l_A(\lambda_i,N)([\zeta])\\ & \in &
\sum^z_{i=1}\Ext^r_A(M, L_i)\Ext^l_A(L_i, N).
\end{array}$$ Conversely, for $[\xi_i] \in \Ext^r_A(M, L_i)$ and
$[\zeta_i] \in \Ext^l_A(L_i, N)$, we have
$$\begin{array}{lll} [\xi_i][\zeta_i] &=& [\xi_i]\Ext^l_A(\rho_i\lambda_i ,
N)([\zeta_i])\\ & =& \Ext^r_A(M,
\lambda_i)([\xi_i])\Ext^l_A(\rho_i,N)([\zeta_i])\\ & \in &
\Ext^r_A(M, \oplus ^z_{i=1}L_i)\Ext^l_A(\oplus ^z_{i=1}L_i, N).
\end{array}$$ \hfill{$\Box$}

Secondly, we have the following ``cancelation lemma'':

\begin{lemma} \label{canlemma} Let $M, N, L_1, L_2, L_3$ be $A$-modules with
$L_2 \cong L_3$. Then, for $r,l \geq 0$, as the $K$-subspaces of
$\Ext^{r+l}_A(M,N)$,
$$\Ext^r_A(M, \oplus ^3_{i=1}L_i)\Ext^l_A(\oplus ^3_{i=1}L_i,
N) = \Ext^r_A(M, \oplus ^2_{i=1}L_i)\Ext^l_A(\oplus ^2_{i=1}L_i,
N).$$ \end{lemma}

{\bf Proof.} Assume that $\sigma : L_2 \rightarrow L_3$ is an
isomorphism of $A$-modules. For $[\xi_3] \in \Ext^r_A(M, L_3)$ and
$[\zeta_3] \in \Ext^l_A(L_3, N)$, we have
$$\begin{array}{lll} [\xi_3][\zeta_3] &=& [\xi_3]\Ext^l_A(\sigma \sigma ^{-1},
N)([\zeta_3])\\ & =& \Ext^r_A(M, \sigma
^{-1})([\xi_3])\Ext^l_A(\sigma ,N)([\zeta_3])\\ & \in & \Ext^r_A(M,
L_2)\Ext^l_A(L_2, N).
\end{array}$$
By Lemma~\ref{sumlemma}, we have
$$\begin{array}{lll}
\Ext^r_A(M, \oplus ^3_{i=1}L_i)\Ext^l_A(\oplus ^3_{i=1}L_i, N)
& = & \sum^3_{i=1}\Ext^r_A(M, L_i)\Ext^l_A(L_i, N)\\
& = & \sum^2_{i=1}\Ext^r_A(M, L_i)\Ext^l_A(L_i, N)\\ & = &
\Ext^r_A(M, \oplus ^2_{i=1}L_i)\Ext^l_A(\oplus ^2_{i=1}L_i, N).
\end{array}$$ \hfill{$\Box$}

For any $K$-algebra $A$, denote by $\Mod A$ (resp. $\mod A$) the
category of (resp. finitely generated) left $A$-modules. We have the
following ``extension group isomorphism lemma'':

\begin{lemma} \label{extlemma} Let $F: \Mod A \rightarrow \Mod A'$
(resp. $F: \mod A \rightarrow \mod A'$) be an equivalence functor.
Then $F$ induces $K$-vector space isomorphisms $F_{MN}:
\Ext^r_A(M,N) \rightarrow \Ext^r_{A'}(F(M),F(N))$ for all $M,N$ in
$\Mod A$ (resp. $\mod A$) and $r \geq 0$.
\end{lemma}

{\bf Proof.} In case $r=0$ the isomorphism is well-known (cf.
\cite[21.2. Proposition]{AF}). By \cite[16.3. Proposition]{AF},
$F(M_1 \oplus M_2)$ is a direct sum $F(M_1)$ and $F(M_2)$ with
injections $F(\lambda_1)$ and $F(\lambda_2)$. Thus there are
$A'$-modules maps $\phi_{M_1M_2} : F(M_1 \oplus M_2) \rightarrow
F(M_1) \oplus F(M_2)$ and $\psi_{M_1M_2} : F(M_1) \oplus F(M_2)
\rightarrow F(M_1 \oplus M_2)$ such that
$\phi_{M_1M_2}\psi_{M_1M_2}=1 , \psi_{M_1M_2}\phi_{M_1M_2} = 1,
\phi_{M_1M_2} F(\lambda_i)=\lambda_i$ and $F(\rho_i)\psi_{M_1M_2} =
\rho_i$ for $i = 1,2$. Hence $\phi_{MM} F(\Delta)=\Delta$ and
$F(\nabla) \psi_{MM} = \nabla$.

Let $f_i : M_i \rightarrow N_i$ be $A$-module maps with $i = 1,2$.
Applying the formulae above, we can show that $\phi_{N_1N_2}
F(diag\{f_1,f_2\})) = diag\{F(f_1),F(f_2)\}\phi_{M_1M_2}$, since
$diag\{f_1,f_2\} = \lambda_1f_1\rho_1 + \lambda_2f_2\rho_2.$

The functor $F$ induces a map $F_{MN}: \Ext^r_A(M,N) \rightarrow
\Ext^r_{A'}(F(M),F(N)),$ \linebreak $[\xi] \rightarrow [F(\xi)]$.
For $[\xi_1], [\xi_2] \in \Ext^r_A(M,N)$, using the preparations
above we can show routinely that $F_{MN}([\xi_1] + [\xi_2]) =
F_{MN}([\xi_1]) + F_{MN}([\xi_2])$. So $F_{MN}$ is a $K$-linear map.

Let $G$ be an inverse equivalence of $F$. Then $G$ also induces a
$K$-linear map $G_{M'N'}: \Ext^r_{A'}(M',N') \rightarrow
\Ext^r_A(G(M'),G(N')),[\zeta] \rightarrow [G(\zeta)]$.

The natural isomorphism $\eta$ from the functor $GF$ to the
identity functor induces an isomorphism $\eta_{MN} :
\Ext^r_A(GF(M),GF(N)) \rightarrow \Ext^r_A(M,N), [\xi] \mapsto
[\eta(\xi)]$, where $\eta(\xi)$ is obtained from $\xi$ by
replacing the monomorphism (resp. epimorphism) in $\xi$ with the
composition of it and the isomorphism $\eta^{-1}_N: N \rightarrow
GF(N)$ (resp. $\eta_M: GF(M) \rightarrow M$).

Since $\eta_{GF(M)GF(N)}G_{F(M)F(N)}F_{MN}=1$, $G_{F(M)F(N)}F_{MN}$
is a bijection. Similarly, $F_{GF(M)GF(N)}G_{F(M)F(N)}$ is a
bijection. Hence $G_{F(M)F(N)}$ is a bijection. It follows that
$F_{MN}$ is a bijection. \hfill{$\Box$}

\begin{theorem} \label{moritaequivalencetheorem} Let
$F: \Mod A \rightarrow \Mod A'$ or $F: \mod A \rightarrow \mod A'$
be an equivalence functor. If $A$ is a quasi-Koszul algebra then so
is $A'$. If $M$ is a quasi-Koszul $A$-module then $F(M)$ is a
quasi-Koszul $A'$-module.
\end{theorem}

{\bf Proof.} By \cite[21.8. Proposition and 27.8. Corollary]{AF},
$A'$ is also noetherian and semiperfect. Denote by $S_1,...,S_s$ a
complete set of the representatives of simple $A$-modules. Then
$S'_1:=F(S_1),...,S'_s:=F(S_s)$ is a complete set of the
representatives of simple $A'$-modules (cf. \cite[\S 21 and \S
27]{AF}). Assume that $A/J \cong \oplus^s_{i=1}S^{u_i}_i=
\oplus^s_{i=1}\oplus^{u_i}_{k=1}S_{ik}$ and $A'/J' \cong
\oplus^s_{j=1}S'^{u'_j}_j= \oplus^s_{j=1}\oplus^{u'_j}_{l=1}S'_{jl}$
where $J'$ is the Jacobson radical of $A'$, $S_{ik}=S_i$ and
$S'_{jl}=S'_j$ for all $i,j,k,l$.

The algebra $A$ is quasi-Koszul, by definition, $\Ext^*_A(A/J,A/J)$
is generated in degree 0 and 1. So is
$\Ext^*_A(\oplus^s_{i=1}S^{u_i}_i,\oplus^s_{i=1}S^{u_i}_i)
=\Ext^*_A(\oplus^s_{i=1}\oplus^{u_i}_{j=1}S_{ij},\oplus^s_{i=1}\oplus^{u_i}_{j=1}S_{ij})$.
Thus $$\begin{array}{rcl}
&&\Ext^{r+1}_A(\oplus^s_{i=1}\oplus^{u_i}_{j=1}S_{ij},
\oplus^s_{i=1}\oplus^{u_i}_{j=1}S_{ij})\\
&=&\Ext^r_A(\oplus^s_{i=1}\oplus^{u_i}_{j=1}S_{ij},
\oplus^s_{i=1}\oplus^{u_i}_{j=1}S_{ij})
\Ext^1_A(\oplus^s_{i=1}\oplus^{u_i}_{j=1}S_{ij},
\oplus^s_{i=1}\oplus^{u_i}_{j=1}S_{ij})
\end{array}$$ for all $r \geq 1$.
By Section~\ref{doeg}, we have
$$\begin{array}{rcl}
&&\coprod^s_{i=1}\coprod^{u_i}_{j=1}\coprod^s_{k=1}\coprod^{u_k}_{l=1}
E^{r+1}_A(S_{ij},S_{kl})\\&=&
\sum^s_{i=1}\sum^{u_i}_{j=1}\sum^s_{k=1}\sum^{u_k}_{l=1}
E^r_A(S_{ij},\oplus^s_{p=1}\oplus^{u_p}_{q=1}S_{pq})
E^1_A(\oplus^s_{p=1}\oplus^{u_p}_{q=1}S_{pq},S_{kl})
\end{array}$$ for all $r \geq 1$.
Note that $E^r_A(S_{ij},\oplus^s_{p=1}\oplus^{u_p}_{q=1}S_{pq})
E^1_A(\oplus^s_{p=1}\oplus^{u_p}_{q=1}S_{pq},S_{kl}) \subseteq
E^{r+1}_A(S_{ij},S_{kl}) $. Thus $$E^{r+1}_A(S_{ij},S_{kl}) =
E^r_A(S_{ij},\oplus^s_{p=1}\oplus^{u_p}_{q=1}S_{pq})
E^1_A(\oplus^s_{p=1}\oplus^{u_p}_{q=1}S_{pq},S_{kl})$$ for all
$i,j,k,l,r$. It follows
$$\Ext^{r+1}_A(S_{ij},S_{kl}) =
\Ext^r_A(S_{ij},\oplus^s_{p=1}\oplus^{u_p}_{q=1}S_{pq})
\Ext^1_A(\oplus^s_{p=1}\oplus^{u_p}_{q=1}S_{pq},S_{kl})$$ for all
$i,j,k,l,r$. By Lemma~\ref{canlemma}, we obtain
$$\Ext^{r+1}_A(S_{ij},S_{kl}) =
\Ext^r_A(S_{ij},\oplus^s_{p=1}S_p)
\Ext^1_A(\oplus^s_{p=1}S_p,S_{kl})$$ for all $i,j,k,l,r$. Applying
Lemma~\ref{extlemma} and Lemma~\ref{canlemma}, we have
$$\begin{array}{lll} &&
\Ext^{r+1}_{A'}(S'_{ij},S'_{kl})\\  & = &
F_{S_{ij}S_{kl}}(\Ext^{r+1}_A(S_{ij},S_{kl}))\\ & = &
F_{S_{ij}S_{kl}}(\Ext^r_A(S_{ij},\oplus^s_{p=1}S_p)
\Ext^1_A(\oplus^s_{p=1}S_p,S_{kl}))\\ & = &
\Ext^r_{A'}(S'_{ij},F(\oplus^s_{p=1}S_p)) \Ext^1_{A'}(F(\oplus^s_{p=1}S_p),S'_{kl}))\\
& = & \Ext^r_{A'}(S'_{ij},\oplus^s_{p=1}S'_p)
\Ext^1_{A'}(\oplus^s_{p=1}S'_p,S'_{kl})\\
& = & \Ext^r_{A'}(S'_{ij},\oplus^s_{p=1}\oplus^{u'_p}_{q=1}S'_{pq})
\Ext^1_{A'}(\oplus^s_{p=1}\oplus^{u'_p}_{q=1}S'_{pq},S'_{kl})
\end{array}$$ for all $i,j,k,l,r$. It follows
$$E^{r+1}_{A'}(S'_{ij},S'_{kl}) =
E^r_{A'}(S'_{ij},\oplus^s_{p=1}\oplus^{u'_p}_{q=1}S'_{pq})
E^1_{A'}(\oplus^s_{p=1}\oplus^{u'_p}_{q=1}S'_{pq},S'_{kl})$$ for all
$i,j,k,l,r$. Thus
$$\begin{array}{rcl}
&&\coprod^s_{i=1}\coprod^{u'_i}_{j=1}\coprod^s_{k=1}\coprod^{u'_k}_{l=1}
E^{r+1}_{A'}(S'_{ij},S'_{kl})\\&=&
\sum^s_{i=1}\sum^{u'_i}_{j=1}\sum^s_{k=1}\sum^{u'_k}_{l=1}
E^r_{A'}(S'_{ij},\oplus^s_{p=1}\oplus^{u'_p}_{q=1}S'_{pq})
E^1_{A'}(\oplus^s_{p=1}\oplus^{u'_p}_{q=1}S'_{pq},S'_{kl})
\end{array}$$ for all $r \geq 1$. It follows from Section~\ref{doeg} that $$\begin{array}{rcl}
&&\Ext^{r+1}_{A'}(\oplus^s_{i=1}\oplus^{u'_i}_{j=1}S'_{ij},\oplus^s_{i=1}\oplus^{u'_i}_{j=1}S'_{ij})\\&=&
\Ext^r_{A'}(\oplus^s_{i=1}\oplus^{u'_i}_{j=1}S'_{ij},\oplus^s_{i=1}\oplus^{u'_i}_{j=1}S'_{ij})
\Ext^1_{A'}(\oplus^s_{i=1}\oplus^{u'_i}_{j=1}S'_{ij},\oplus^s_{i=1}\oplus^{u'_i}_{j=1}S'_{ij})
\end{array}$$ for all $r \geq
1$. So
$\Ext^*_{A'}(\oplus^s_{i=1}S'^{u'_i}_i,\oplus^s_{i=1}S'^{u'_i}_i)=
\Ext^*_{A'}(\oplus^s_{i=1}\oplus^{u'_i}_{j=1}S'_{ij},\oplus^s_{i=1}\oplus^{u'_i}_{j=1}S'_{ij})$
is generated in degree 0 and 1. Thus $\Ext^*_{A'}(A'/J',A'/J')$ is
generated in degree 0 and 1, i.e., $A'$ is quasi-Koszul.

\medskip

Now suppose $M$ is a quasi-Koszul $A$-module, we show that $F(M)$ is
a quasi-Koszul $A'$-module. By definition,
$$\Ext^{r+1}_A(M,A/J)=\Ext^r_A(M,A/J)\Ext^1_A(A/J,A/J)$$ for all $r
\geq 0$. It follows from Lemma~\ref{extlemma} that
$$\Ext^{r+1}_{A'}(F(M),F(A/J))=\Ext^r_{A'}(F(M),F(A/J))\Ext^1_{A'}(F(A/J),F(A/J))$$
for all $r \geq 0$. Therefore
$$\begin{array}{rcl}
&&\Ext^{r+1}_{A'}(F(M),\oplus^s_{i=1}\oplus^{u_i}_{j=1}S'_{ij})\\
&=&\Ext^r_{A'}(F(M),\oplus^s_{i=1}\oplus^{u_i}_{j=1}S'_{ij})
\Ext^1_{A'}(\oplus^s_{i=1}\oplus^{u_i}_{j=1}S'_{ij},\oplus^s_{i=1}\oplus^{u_i}_{j=1}S'_{ij})
\end{array}$$ for all $r \geq 0$.
By Section~\ref{doeg}, we have
$$\begin{array}{rcl}
&&\coprod^s_{i=1}\coprod^{u_i}_{j=1} E^{r+1}_{A'}(F(M),S'_{ij})\\&=&
\sum^s_{i=1}\sum^{u_i}_{j=1}
\Ext^r_{A'}(F(M),\oplus^s_{p=1}\oplus^{u_p}_{q=1}S'_{pq})
E^1_{A'}(\oplus^s_{p=1}\oplus^{u_p}_{q=1}S'_{pq},S'_{ij})
\end{array}$$ for all $r \geq 0$. Note that
$\Ext^r_{A'}(F(M),\oplus^s_{p=1}\oplus^{u_p}_{q=1}S'_{pq})
E^1_{A'}(\oplus^s_{p=1}\oplus^{u_p}_{q=1}S'_{pq},S'_{ij}) \subseteq$
\linebreak $E^{r+1}_{A'}(F(M),S'_{ij})$. Thus
$$E^{r+1}_{A'}(F(M),S'_{ij}) =
\Ext^r_{A'}(F(M),\oplus^s_{p=1}\oplus^{u_p}_{q=1}S'_{pq})
E^1_{A'}(\oplus^s_{p=1}\oplus^{u_p}_{q=1}S'_{pq},S'_{ij})$$ for all
$i,j,r$. It follows $$\Ext^{r+1}_{A'}(F(M),S'_{ij}) =
\Ext^r_{A'}(F(M),\oplus^s_{p=1}\oplus^{u_p}_{q=1}S'_{pq})
\Ext^1_{A'}(\oplus^s_{p=1}\oplus^{u_p}_{q=1}S'_{pq},S'_{ij})$$ for
all $i,j,r$. By Lemma~\ref{canlemma}, we have
$$\begin{array}{rcl} &&\Ext^{r+1}_{A'}(F(M),S'_{ij})\\ & = &
\Ext^r_{A'}(F(M),\oplus^s_{p=1}S'_p)
\Ext^1_{A'}(\oplus^s_{p=1}S'_p,S'_i)\\ & = &
\Ext^r_{A'}(F(M),\oplus^s_{p=1}\oplus^{u'_p}_{q=1}S'_{pq})
\Ext^1_{A'}(\oplus^s_{p=1}\oplus^{u'_p}_{q=1}S'_{pq},S'_{ij})
\end{array}$$ for all $i,j,r$. It follows $$E^{r+1}_{A'}(F(M),S'_{ij}) =
\Ext^r_{A'}(F(M),\oplus^s_{p=1}\oplus^{u'_p}_{q=1}S'_{pq})
E^1_{A'}(\oplus^s_{p=1}\oplus^{u'_p}_{q=1}S'_{pq},S'_{ij})$$ for all
$i,j,r$. Thus
$$\begin{array}{ll}
&\coprod^s_{i=1}\coprod^{u'_i}_{j=1} E^{r+1}_{A'}(F(M),S'_{ij})\\
=& \sum^s_{i=1}\sum^{u'_i}_{j=1}
\Ext^r_{A'}(F(M),\oplus^s_{p=1}\oplus^{u'_p}_{q=1}S'_{pq})
E^1_{A'}(\oplus^s_{p=1}\oplus^{u'_p}_{q=1}S'_{pq},S'_{ij})
\end{array}$$ for all $r \geq 0$. By Section~\ref{doeg}, we obtain $$\begin{array}{ll}
&\Ext^{r+1}_{A'}(F(M),\oplus^s_{i=1}\oplus^{u'_i}_{j=1}S'_{ij})\\=&
\Ext^r_{A'}(F(M),\oplus^s_{i=1}\oplus^{u'_i}_{j=1}S'_{ij})
\Ext^1_{A'}(\oplus^s_{i=1}\oplus^{u'_i}_{j=1}S'_{ij},\oplus^s_{i=1}\oplus^{u'_i}_{j=1}S'_{ij})
\end{array}$$ for all $r \geq0$. Thus $$\Ext^{r+1}_{A'}(F(M),A'/J') =
\Ext^r_{A'}(F(M),A'/J') \Ext^1_{A'}(A'/J',A'/J')$$ for all $r \geq
0$, i.e., $F(M)$ is a quasi-Koszul $A'$-module. \hfill{$\Box$}

\begin{remark} Though the quasi-Koszulity of algebras is invariant
under Morita equivalence, it is not so under derived equivalence
even tilting equivalence in general. Indeed, all path algebras are
Koszul thus quasi-Koszul. If all their tilted algebras were
quasi-Koszul then, as length graded algebras, all monomial tilted
algebras would be Koszul, thus quadratic \cite[Corollary
7.3]{GM1996}. It is a contradiction \cite[Appendix]{R1099}.
Nevertheless, since the global dimension of tilted algebras are less
equal to 2 \cite[Theorem (5.2)]{HR}, all quadratic tilted algebras
are Koszul \cite[Theorem 7.2]{GM1996}.
\end{remark}

\begin{corollary} \label{matrixcorollary} $A$ is quasi-Koszul
if and only if so is $M_n(A)$. \end{corollary}

\section{Skew group algebras and quasi-Koszulity} \label{skewgroupalgebra}

In this section, we show that, a finite-dimensional $K$-algebra $A$
with an action of $G$ is quasi-Koszul if and only if so is the skew
group algebra $A \ast G$, where $G$ is a finite group satisfying
$\char K \nmid |G|$. Moreover, we prove that the Yoneda algebra
$\Ext^*_{A \ast G}((A \ast G)/J(A \ast G),(A \ast G)/J(A \ast G))$
of the skew group algebra $A \ast G$ is isomorphic to the skew group
algebra $\Ext^*_A(A/J,A/J) \ast G$ of the Yoneda algebra
$\Ext^*_A(A/J,A/J)$. The results of this section are essentially due
to Mart\'{i}nez-Villa (cf. \cite{Ms}).

\subsection{Skew group algebras}

Let $A$ be a $K$-algebra and $G$ a group of $K$-algebra
automorphisms of $A$. Then the {\it skew group algebra} $A \ast G$
is defined as follows: As a $K$-vector space $A \ast G = A \otimes
_K KG$. For $a \in A$ and $g \in G$, write $ag$ instead of $a
\otimes g$, and define the multiplication by $ga:=g(a)g$.

Let $A$ be a $K$-algebra and $G$ a group. Then an {\it action} of
$G$ on $A$ is a map $G \times A \rightarrow A, (g,a) \mapsto ga$
satisfying $(g_1g_2)a=g_1(g_2a), 1a=a, g(a_1+a_2)=ga_1+ga_2,
g(a_1a_2)=(ga_1)(ga_2), g(ka)=k(ga),$ for all $g, g_1,g_2 \in G,
a,a_1,$ \linebreak $a_2 \in A, k \in K$.

For a $K$-algebra $A$, giving a group $G$ of $K$-automorphisms of
$A$ is equivalent to giving an action of a group $G$ on $A$, or
giving a group homomorphism $G \rightarrow \Aut_k(A)$ where
$\Aut_k(A)$ denotes the group of $K$-algebra automorphisms of $A$,
or giving a $KG$-module algebra structure (cf. \cite[Proposition
1.2]{CM}).

\subsection{Quasi-Koszulity of skew group algebras}

Let $G$ be a finite group and $M$ a $KG$-module. Denote by $M^G$
the $KG$-module $\{m \in M | gm=m \mbox{ for all } g \in G\}$. If
$\char K \nmid |G|$ then the fixed point functor $(-)^G: \Mod KG
\rightarrow \Mod KG$ is exact (cf. \cite[Lemma 3]{Ms}).

\begin{lemma} (cf. \cite[Lemma 4]{Ms}) \label{actlemma}
Let $A$ be a $K$-algebra, $G$ a finite group acting on $A$, and
$M,N,L$ three $A \ast G$-modules. Then

(1) $\Hom_A(M,N)$ is a $KG$-module defined by $(g
\phi)(m)=g\phi(g^{-1}m)$ for $ g \in G, \phi \in \Hom_A(M,N), m
\in M$.

(2) $\Hom_A(M,N)^G=\Hom_{A \ast G}(M,N)$.

(3) $\Ext^i_A(M,N)$ is a $KG$-module satisfying
$g([\xi][\zeta])=(g[\xi])(g[\zeta])$ for $g \in G$, $[\xi] \in
\Ext^j_A(L,M)$ and $[\zeta] \in \Ext^i_A(M,N)$.
\end{lemma}

{\bf Proof.} (1), (2): Easy to check.

(3): To each $g \in G$, we associate a functor $(-)^g : \Mod A
\rightarrow \Mod A$. For each $X$ in $\Mod A$, $X^g$ is defined as
follows: As a $K$-vector space $X^g=X$. For $a \in A$ and $x \in
X^g$, $a \cdot x :=(ga)x$. For $X,Y$ in $\Mod A$ and $\phi \in
\Hom_A(X,Y)$, $\phi^g := \phi$. Obviously, the functor $(-)^g$ is
not only an exact functor but also an automorphism of $\Mod A$ with
the inverse $()^{g^{-1}}$.

Since $M$ is an $A \ast G$-module, we have an $A$-module
isomorphism $\Psi^g_M : M \rightarrow M^g, m \mapsto gm$ with the
inverse $\Psi^{g^{-1}}_{M^g} : M^g \rightarrow M$ (cf.
\cite[Proposition 2.5]{MMM}). Note that for $\psi \in \Hom_A(M,N)$
we have $\Psi^g_{N^{g^{-1}}} \psi^{g^{-1}} \Psi^{g^{-1}}_M = g
\psi$.

For $g \in G$ and $[\xi] \in \Ext^i_A(M,N)$, the exact functor
$(-)^{g^{-1}}$ provides an $i$-extension $\xi^{g^{-1}}$ of
$M^{g^{-1}}$ by $N^{g^{-1}}$. Since $M$ and $N$ are $A \ast
G$-modules, replacing the epimorphism (resp. monomorphism) in
$\xi^{g^{-1}}$ with the composition of it and the $A$-module
isomorphism $\Psi^g_{M^{g^{-1}}}$ (resp. $\Psi^{g^{-1}}_N$), we
obtain an $i$-extension $g\xi$ of $M$ by $N$. One can shows that
$\Ext^i_A(M,N)$ is a $KG$-module defined by $g [\xi]:=[g\xi]$ and
satisfying $g([\xi][\zeta])=(g[\xi])(g[\zeta])$ for $g \in G$,
$[\xi] \in \Ext^j_A(L,M)$ and $[\zeta] \in \Ext^i_A(M,N)$.
\hfill{$\Box$}

\begin{lemma} (cf. \cite[Corollary 5]{Ms}) \label{leslemma}
Let $A$ be a $K$-algebra, $G$ a finite group acting on $A$, and
$\varsigma: 0 \rightarrow L \stackrel{\phi}{\rightarrow} M
\stackrel{\psi}{\rightarrow} N \rightarrow 0$ an exact sequence of
$A \ast G$-modules. Then for any $A \ast G$-module $X$ two long
exact sequences
$$\begin{array}{lll} 0 & \rightarrow & \Hom_A(X,L) \rightarrow
\Hom_A(X,M) \rightarrow \Hom_A(X,N)\\ & \rightarrow & \Ext^1_A(X,L)
\rightarrow \Ext^1_A(X,M) \rightarrow \Ext^1_A(X,N) \rightarrow
\cdots \\ &&\\ 0 & \rightarrow & \Hom_A(N,X) \rightarrow \Hom_A(M,X)
\rightarrow \Hom_A(L,X)\\ & \rightarrow & \Ext^1_A(N,X) \rightarrow
\Ext^1_A(M,X) \rightarrow \Ext^1_A(L,X) \rightarrow \cdots
\end{array}$$ are exact sequences of $KG$-modules and $KG$-module
maps. \end{lemma}

{\bf Proof.} By Lemma~\ref{actlemma} (1) and (3), we know that all
terms in the long exact sequences are $KG$-modules.

Since $\phi$ is an $A \ast G$-module map, by Lemma~\ref{actlemma}
(2), we have $g \phi = \phi$ for all $g \in G$. For any $[\xi] \in
\Ext^i_A(X,L)$, using $\Psi^g_{M^{g^{-1}}} \phi^{g^{-1}}
\Psi^{g^{-1}}_L = g \phi$, we can show that $g \Ext^i_A(X,
\phi)([\xi]) = \Ext^i_A(X, g \phi)([g\xi]) = \Ext^i_A(X,
\phi)(g[\xi])$. It follows that $\Ext^i_A(X, \phi)$ is a $KG$-module
map. Similar for $\Ext^i_A(X, \psi)$.

Let $\delta : \Ext^i_A(X,N) \rightarrow \Ext^{i+1}_A(X,L)$ be the
connecting map and $[\zeta] \in \Ext^i_A(X,N)$. Since $\varsigma$
is an exact sequence of $A \ast G$-modules, we have $g [\varsigma]
= [\varsigma]$. By Lemma~\ref{actlemma} (3), we have $\delta (g
[\zeta]) = (g [\zeta])[\varsigma] = (g [\zeta]) (g [\varsigma]) =
g ([\zeta \varsigma]) = g \delta([\zeta])$. It follows that
$\delta$ is a $KG$-module map.

Hence the first exact sequence is an exact sequence of
$KG$-modules and $KG$-module maps. Similar for the second exact
sequence. \hfill{$\Box$}

\begin{lemma} (cf. \cite[Lemma 6]{Ms}) \label{projlemma}
Let $A$ be a $K$-algebra, $G$ a finite group acting on $A$ with
$\char K \nmid |G|$ and $P$ a finitely generated $A \ast
G$-module. Then $P$ is projective if and only if $P$ is projective
as an $A$-module. \end{lemma}

{\bf Proof.} If $P$ is projective as an $A \ast G$-module then $P$
is a direct summand of a free $A \ast G$-module $(A \ast G)^i$. So
$P$ is a direct summand of a free $A$-module $(A \ast G)^i \cong
A^{|G|i}$. Hence $P$ is a projective $A$-module.

Conversely, if $P$ is a projective $A$-module then by
Lemma~\ref{leslemma} the functor $\Hom_A(P,-) : \Mod A \ast G
\rightarrow \Mod KG$ is exact. Since the functor $(-)^G : \Mod KG
\rightarrow \Mod KG$ is exact, the functor $\Hom_A(P,-)^G : \Mod A
\ast G \rightarrow \Mod KG$ is also exact. By Lemma~\ref{actlemma}
(2) the functor $\Hom_{A \ast G}(P,-)=$ \linebreak $\Hom_A(P,-)^G
: \Mod A \ast G \rightarrow \Mod KG$ is exact. Thus $P$ is also
projective as an $A \ast G$-module. \hfill{$\Box$}

\medskip

Let $M$ be an $A \ast G$-module and $W$ a $KG$-module. Then $M
\otimes_K W$ is an $A \ast G$-module defined by $(ag)(m \otimes
w)= agm \otimes gw$ for $ a \in A, g \in G, m \in M, w \in W$.

\begin{lemma} (cf. \cite[Lemma 8]{Ms}) \label{homlemma}
Let $A$ be a $K$-algebra, $G$ a finite group acting on $A$, $P$ a
finitely generated projective $A \ast G$-module, $N$ an $A \ast
G$-module and $W$ a $KG$-module. Then we have a $K$-vector space
isomorphism $\theta_P : \Hom_A(P,N) \otimes _K W \rightarrow
\Hom_{A \ast G}(P \otimes _K KG, N \otimes _K W)$ defined by
$\theta_P(\phi \otimes w)(p \otimes g):=(g \phi)(p) \otimes gw$.
\end{lemma}

{\bf Proof.} Consider the natural isomorphisms $\sigma :
\Hom_A(A,N) \otimes_K W \rightarrow N \otimes _K W, \phi \otimes w
\mapsto \phi(1) \otimes w$ and $\tau : \Hom_{A \ast G}(A \otimes
_K KG, N \otimes _K W) \rightarrow N \otimes _K W, \psi \mapsto
\psi(1 \otimes 1)$. We have $\theta_A=\tau ^{-1} \sigma$. Thus
Lemma~\ref{homlemma} holds for $P=A$. Next we can prove routinely
that Lemma~\ref{homlemma} holds for all $A^n$ with $n \geq 1$ and
all finitely generated projective $A \ast G$-module $P$.
\hfill{$\Box$}

\begin{lemma} (cf. \cite[Proposition 9]{Ms}) \label{extisolemma}
Let $A$ be a $K$-algebra, $G$ a finite group acting on $A$ with
$\char K \nmid |G|$, $M$ an $A \ast G$-module admitting a finitely
generated $A \ast G$-module projective resolution $\mathbb{P} =
(P_i, \psi_i)_{i \geq 0}$, $N$ an $A \ast G$-module and $W$ a
$KG$-module. Then we have $K$-vector space isomorphisms
$\bar{\theta}_i : \Ext^i_A(M,N) \otimes _K W \rightarrow \Ext^i_{A
\ast G}(M \otimes _K KG, N \otimes _K W)$ defined by
$\bar{\theta}_i(\bar{\phi} \otimes w):=\overline{\theta_{P_i}(\phi
\otimes w)}$, for all $i \geq 0$. \end{lemma}

\begin{remark} \label{correspondenceremark}
Here we view the elements in $\Ext^i_A(M,N)$ as the residue
classes in $\Hom_A(P_i,N)/\Im \Hom_A(\psi_i,N)$ or
$\Hom_A(\Omega^i(M),N)/\Im \Hom_A(\iota_i,N)$ where $\Omega^i(M)$
denotes the $i$-th syzygy of $M$ and $\iota_i : \Omega^i(M)
\hookrightarrow P_{i-1}$ is the natural embedding (cf.
\cite[Chapter III, Theorem 6.4]{Mac}), because in this situation
this viewpoint is more intuitive.
\end{remark}

{\bf Proof of Lemma~\ref{extisolemma}.} By Lemma~\ref{projlemma},
$\mathbb{P}$ is also a finitely generated $A$-module projective
resolution of $M$. Thus $\mathbb{P} \otimes _K KG$ is an exact
sequence of $A \ast G$-modules where the $A \ast G$-module structure
of $P_i \otimes_KKG$ is given by $(ag)(p_i \otimes h)=agp_i \otimes
gh$. By Lemma~\ref{projlemma} again, $\mathbb{P} \otimes _K KG=(P_i
\otimes _K KG, \psi_i \otimes 1)_{i \geq 0}$ is a finitely generated
$A \ast G$-module projective resolution of $M \otimes _K KG$. It
follows from Lemma~\ref{homlemma} that $(\theta_{P_i})_{i \geq 0}$
is a chain isomorphism from the complex $\Hom_A(\mathbb{P},N)
\otimes _K W$ to the complex $\Hom_A(\mathbb{P} \otimes _K KG, N
\otimes _K W)$. This chain isomorphism induces isomorphisms
$\bar{\theta}_i : \Ext^i_A(M,N) \otimes _K W \rightarrow \Ext^i_{A
\ast G}(M \otimes _K KG, N \otimes _K W)$ between the homologies of
these complexes defined by $\bar{\theta}_i(\bar{\phi} \otimes
w):=\overline{\theta_{P_i}(\phi \otimes w)}$ for all $i \geq 1$. The
case of $i=0$ follows from Five Lemma. \hfill{$\Box$}

\begin{theorem} (cf. \cite[Theorem 10]{Ms}) \label{skewtheorem}
Let $A$ be a finite-dimensional $K$-algebra, $G$ a finite group
acting on $A$ with $\char K \nmid |G|$ and $M$ a finitely generated
$A \ast G$-module. Then $M \otimes _KKG$ is a quasi-Koszul $A \ast
G$-module if and only if $M$ is a quasi-Koszul $A$-module. In
particular, the algebra $A \ast G$ is quasi-Koszul if and only if so
is $A$. Moreover, $\Ext^*_{A \ast G}((A \ast G)/J(A \ast G), (A \ast
G)/J(A \ast G)) \cong \Ext^*_A(A/J,A/J) \ast G$ as positively graded
$K$-algebras.
\end{theorem}

{\bf Proof.} Let $\mathbb{P} = (P_k, \psi_k)_{k \geq 0}$ be a
finitely generated $A \ast G$-module projective resolution of $M$.
First of all, we show that the following diagram is commutative for
all $i,j \geq 0$: {\tiny
$$\begin{array}{ccc} (\Ext^i_A(M,A/J) \otimes_KKG)
\times (\Ext^j_A(A/J,A/J) \otimes_KKG) &
\stackrel{\mu}{\rightarrow} & \Ext^{i+j}_A(M,A/J) \otimes_KKG\\
\downarrow
(\bar{\theta}_i,\bar{\theta}_j) && \downarrow \bar{\theta}_{i+j}\\
\Ext^i_{A \ast G}(M \otimes_KKG,(A/J) \ast G) \times \Ext^j_{A \ast
G}((A/J) \ast G,(A/J) \ast G) & \stackrel{\nu}{\rightarrow} &
\Ext^{i+j}_{A \ast G}(M \otimes_KKG,(A/J) \ast G)
\end{array}$$} where $\nu$ is Yoneda product and $\mu$ is defined
by $\mu(\bar{\az} \otimes g, \bar{\bz} \otimes
h):=\bar{\az}(g\bar{\bz}) \otimes gh$.

Let $\mathbb{Q}=(Q_k, \phi_k)_{k \geq 0}$ be a finitely generated $A
\ast G$-module projective resolution of $A/J$. Let $\bar{\az} \in
\Ext^i_A(M,A/J)$ and $\bar{\bz} \in \Ext^j_A(A/J,A/J)$ with $\az \in
\Hom_A(\Omega^i(M),A/J)$ or $\Hom_A(P_i,A/J)$ and $\bz \in
\Hom_A(\Omega^j(A/J),A/J)$ or $\Hom_A(Q_j,A/J)$ (cf.
Remark~\ref{correspondenceremark}). Then $\bar{\az}\bar{\bz}=
\overline{\bz \Omega^j(\az)}$.

Since $\mathbb{P} = (P_k, \psi_k)_{k \geq 0}$ is an $A \ast
G$-module projective resolution of $M$, $\mathbb{P}(i):= (P_{k+i},
\psi_{k+i})_{k \geq 0}$ is an $A \ast G$-module projective
resolution of $\Omega^i(M)$ for all $i \geq 0$. Furthermore,
$\mathbb{P}(i) \otimes _KKG$ is an $A \ast G$-module projective
resolution of $\Omega^i(M) \otimes _KKG$ for all $i \geq 0$. Suppose
the map $\az : \Omega^i(M) \rightarrow A/J$ induces a chain map
$(\az_j)_{j \geq 0}: \mathbb{P}(i) \rightarrow \mathbb{Q}$. Then
$(\theta_{P_{j+i}}(\az_j \otimes g))_{j \geq 0}: \mathbb{P}(i)
\otimes _KKG \rightarrow \mathbb{Q} \otimes _KKG$ is exactly a chain
map induced by the map $\theta_{P_i}(\az \otimes g): P_i \otimes_KKG
\rightarrow (A/J) \otimes_KKG$. It follows that $\Omega
^j(\theta_{P_i}(\az \otimes g))= \theta_{P_{j+i}}(\Omega^j(\az)
\otimes g)$. Thus for $x \otimes l \in P_{j+i} \otimes _KKG$
$$\begin{array}{lll} &&\theta_{Q_j}(\bz \otimes h)\Omega^j(\theta_{P_i}(\az
\otimes g))(x \otimes l)\\ & = & \theta_{Q_j}(\bz \otimes
h)\theta_{P_{j+i}}(\Omega^j(\az) \otimes g))(x \otimes l)\\ & = &
\theta_{Q_j}(\bz \otimes h)((l \Omega^j(\az))(x) \otimes lg)\\
& = & (lg \bz)((l \Omega^j(\az))(x)) \otimes lgh\\ & = & lg \bz
((lg)^{-1}l (\Omega^j(\az)(l^{-1}x))) \otimes lgh\\ & = & l(g \bz
(g^{-1}( \Omega^j(\az)(l^{-1}x)))) \otimes lgh\\ & = & (l((g \bz)
\Omega^j(\az)))(x) \otimes lgh\\ & = & \theta_{P_{j+i}} ((g \bz)
\Omega^j(\az) \otimes gh)(x \otimes l),
\end{array}$$
 i.e., $\theta_{Q_j}(\bz
\otimes h)\Omega^j(\theta_{P_i}(\az \otimes g)) = \theta
_{P_{j+i}}((g \bz) \Omega^j(\az) \otimes gh)$.

Therefore $\bar{\theta}_{i+j}((\bar{\az} \otimes g)(\bar{\bz}
\otimes h)) = \bar{\theta}_{i+j}(\bar{\az}(g\bar{\bz}) \otimes gh) =
\bar{\theta}_{i+j}(\overline{(g \bz) \Omega ^j (\az)} \otimes gh) =
\overline{\theta_{P_{j+i}} ((g \bz) \Omega^j(\az) \otimes gh)} =
\overline{\theta_{Q_j}(\bz \otimes h)\Omega^j(\theta_{P_i}(\az
\otimes g))} = \overline{\theta_{P_i}(\az \otimes g)}$ \linebreak
$\overline{\theta_{Q_j}(\bz \otimes h)} = \bar{\theta}_i(\bar{\az}
\otimes g) \bar{\theta}_j(\bar{\bz} \otimes h)$. So the diagram
above is commutative.

If $\char K \nmid |G|$ then, by \cite[Theorem 1.1]{RR}, $J \ast G$
is the Jacobson radical of $A \ast G$. Moreover, $(A \ast G)/(J \ast
G) \cong (A /J) \ast G$. Thus $M \otimes_KKG$ is a quasi-Koszul $A
\ast G$-module if and only if $$\begin{array}{lll} & & \Ext^{i+j}_{A
\ast G}(M \otimes_KKG,(A/J) \ast G)\\ & = & \Ext^i_{A \ast G}(M
\otimes_KKG,(A/J) \ast G)\Ext^j_{A \ast G}((A/J) \ast G,(A/J) \ast
G) \end{array}$$ for all $i,j \geq 0$, if and only if
$$\begin{array}{lll} & & \Ext^{i+j}_A (M ,A/J) \otimes _KKG\\ & = &
(\Ext^i_A (M ,A/J) \otimes _KKG)(\Ext^j_A (A/J,A/J) \otimes _KKG)\\
& = & (\Ext^i_A (M,A/J) \Ext^j_A (A/J,A/J)) \otimes _KKG
\end{array}$$ for all $i, j \geq 0$, if and only if $$\Ext^{i+j}_A
(M,A/J) = \Ext^i_A (M,A/J) \Ext^j_A (A/J,A/J)$$ for all $i,j \geq
0$, if and only if $M$ is a quasi-Koszul $A$-module. In particular,
taking $M=A/J$ we have the algebra $A \ast G$ is quasi-Koszul if and
only if so is $A$. Moreover, according to the commutative diagram
above, the map $(\bar{\theta}_i)_{i \geq 0}$ is a positively graded
$K$-algebra isomorphism $\Ext^*_A(A/J,A/J) \ast G \cong \Ext^*_{A
\ast G}((A \ast G)/J(A \ast G), (A \ast G)/J(A \ast G))$.
\hfill{$\Box$}

\section{Smash product and quasi-Koszulity} \label{smashproduct}

In this section we show in two ways that a finite-dimensional
$G$-graded $K$-algebra $A$ is quasi-Koszul if and only if so is the
smash product $A \# G^*$, where $G$ is a finite group satisfying
$\char K \nmid |G|$. Moreover, we prove that the Yoneda algebra
$\Ext^*_{A \# G^*}((A \# G^*)/J(A \# G^*), (A \# G^*)/J(A \# G^*))$
of the smash product $A \# G^*$ is isomorphic to the smash product
$\Ext^*_A(A/J, A/J) \# G^*$ of the Yoneda algebra
$\Ext^*_A(A/J,A/J)$ as positively graded $K$-algebras.

\subsection{Smash product}

Let $G$ be a finite group and $A= \coprod_{g \in G} A_g$ be a
$G$-graded $K$-algebra. A {\it graded $A$-module} $M$ is an
$A$-module, together with a direct sum decomposition
$M=\coprod_{g\in G}M_g$ of $K$-vector spaces such that $A_gM_h
\subseteq M_{gh}$ for all $g, h \in G$. Denote by $\Gr A$ the
category whose objects are all graded $A$-modules and whose
morphisms $\phi: M \rightarrow N$ are morphisms in $\Mod A$ such
that $\phi(M_g) \subset N_g$ for all $g\in G$. The category $\Gr A$
is not only an abelian category but also a Grothendieck category
(cf. \cite[Section 2.2]{NV}). Since it is closed under kernel and
has enough projective objects, each object in $\Gr A$ has a graded
projective resolution.

Let $KG^*$ be the dual algebra of $KG$ and $\{p_g|g \in G\}$ its
dual basis, i.e., for $g \in G$ and $x=\sum_{h \in G}a_hh \in KG$
one has $p_g(x)=a_g \in K$ and $p_gp_h=\delta_{gh}p_h$ where
$\delta_{gh}$ is the Kronecker delta. The {\it smash product},
denoted by $A \# G^*$ (cf. \cite[Section 1]{CM} and \cite[Section
2]{B}), is the $K$-vector space $A \otimes _K KG^*$ with the
multiplication given by $(a \# p_g)(b \# p_h) := ab_{gh^{-1}} \#
p_h$ where $a \# p_g$ denotes the element $a \otimes p_g$.

For graded $A$-modules $M$ and $N$, we denote by $\HOM_A(M, N)$ the
set $\oplus_{g\in G}\Hom_{A}(M, N)_g$, where $\Hom_{A}(M, N)_g=\{f
\in \Hom_A(M, N)| f(M_h) \subset N_{hg}$ for all $h \in G\}$. We
write $\EXT^i_A(M, N)$ for the derived functor of $\HOM_A(M, N)$.
Since $G$ is finite, we have $\HOM_A(M, N)=\Hom_A(M, N)$ and
$\EXT^i_A(M, N)=\Ext^i_A(M, N)$ (cf. \cite[Corollary 2.4.6 and
Corollary 2.4.7]{NV}). Now that $\Ext^i_A(M, N)$ may be computed
from a graded projective resolution of $M$, the grading on each
$\Hom_A(M, N)$ induces a grading on $\Ext^i_A(M, N)$. Therefore
$\Ext_A^*(A/J, A/J) \# G^*$ is well-defined.

For any graded $A$-module $M$, $M \otimes _KKG^*$ is a left $A \#
G^*$-module defined by $(a \otimes p_g)(m \otimes p_h)=am_{gh^{-1}}
\otimes p_h$ for all $a \in A, m \in M$ and $g,h \in G$. If $P$ is a
finitely generated graded projective $A$-module then $P \otimes
_KKG^*$ is a finitely generated projective $A \# G^*$-module.

\subsection{Quasi-Koszulity of smash product}

\begin{lemma} \label{homsmashlemma} Let $G$ be a finite group,
$A$ a $G$-graded $K$-algebra, $P$ a finitely generated graded
projective $A$-module and $N$ a graded $A$-module. Then there is a
$K$-vector space isomorphism $\theta_P: \Hom_A(P,N) \otimes_KKG^*
\rightarrow \Hom_{A \# G^*}(P \otimes_KKG^*, N\otimes_KKG^*)$
defined by $\theta_P(\phi \otimes p_g)(x_h \otimes
p_l)=\phi(x_h)_{hlg^{-1}}\otimes\,p_g$ for all $x_h \in P_h$ and
$p_g, p_l \in KG^*$.
\end{lemma}

{\bf Proof.} Consider the natural isomorphism $\sigma :
\Hom_A(A,N) \otimes_KKG^* \rightarrow N \otimes_KKG^*, \phi
\otimes p_g \mapsto \phi(1) \otimes p_g$ and $\tau : \Hom(A
\otimes_KKG^*, N \otimes_KKG^*) \rightarrow N \otimes_KKG^*, \psi
\mapsto \psi(1 \otimes 1)$. We have $\theta_A=\tau ^{-1} \sigma$.
Thus Lemma~\ref{homsmashlemma} holds for $P=A$. Next we can prove
routinely that Lemma~\ref{homsmashlemma} holds for all $A^n$ with
$n \geq 1$ and all finitely generated projective $A \ast G$-module
$P$. \hfill{$\Box$}

\begin{lemma} \label{extisosmashlemma} Let $G$ be a finite group,
$A$ a $G$-graded $K$-algebra, $M$ a graded $A$-module admitting a
finitely generated graded projective resolution $\mathbb{P}=(P_i,
\psi_i)_{i \geq 0}$ and $N$ a graded $A$-module. Then there is a
natural isomorphism $\overline{\theta}_i: \Ext_A^i(M, N)
\otimes_KKG^* \rightarrow \Ext^i_{A \# G^*}(M \otimes_KKG^*, N
\otimes_KKG^*)$ defined by
$\overline{\theta}_i(\overline{\phi}\otimes
p_g):=\overline{\theta_{P_i}(\phi \otimes p_g)}$ for each $i \geq
0$.
\end{lemma}

{\bf Proof.} Note that $\mathbb{P} \otimes_KKG^*=(P_i \otimes_KKG^*,
\psi_i \otimes 1)_{i \ge 0}$ is a finitely generated $A \#
G^*$-module projective resolution of $M \otimes_KKG^*$. It follows
from Lemma~\ref{homsmashlemma} that $(\theta_{P_i})_{i \geq 0}$ is a
chain isomorphism from the complex $\Hom_A(\mathbb{P}, N)
\otimes_KKG^*$ to the complex $\Hom_{A \#
G^*}(\mathbb{P}\otimes_KKG^*, N\otimes_KKG^*)$. This chain
isomorphism induces isomorphisms $\overline{\theta}_i: \Ext_A^i(M,
N) \otimes_KKG^* \rightarrow \Ext^i_{A \# G^*}(M\otimes_KKG^*, N
\otimes_KKG^*)$ between the homologies of these two complexes
defined by $\overline{\theta}_i(\overline{\phi}\otimes
p_g):=\overline{\theta_{P_i}(\phi\otimes p_g)}$ for all $i \geq 1$.
The case $i=0$ follows from Five Lemma. \hfill{$\Box$}

\begin{theorem} \label{smashtheorem}
Let $A$ be a finite-dimensional $G$-graded $K$-algebra, where $G$ is
a finite group satisfying $\char K \nmid |G|$. Let $M$ be a finitely
generated graded $A$-module. Then $M \otimes_KKG^*$ is a
quasi-Koszul $A\#G^*$-module if and only if $M$ is a quasi-Koszul
$A$-module. In particular, $A \# G^*$ is quasi-Koszul if and only if
so is $A$. Moreover, $\Ext^*_{A \# G^*}((A \# G^*)/J(A \# G^*), (A
\# G^*)/J(A \# G^*)) \cong \Ext^*_A(A/J, A/J) \# G^*$ as positively
graded $K$-algebras.
\end{theorem}

{\bf Proof.} Let $\mathbb{P}=(P_k, \psi_k)_{k \ge 0}$ be a finitely
generated graded projective resolution of $M$. First of all, we show
that the following diagram is commutative: {\tiny
$$\begin{array}{ccc} (\Ext^i_A(M,A/J) \otimes_KKG^*)
\times (\Ext^j_A(A/J,A/J) \otimes_KKG^*) &
\stackrel{\mu}{\rightarrow} & \Ext^{i+j}_A(M,A/J) \otimes_KKG^*\\
\downarrow
(\bar{\theta}_i,\bar{\theta}_j) && \downarrow \bar{\theta}_{i+j}\\
\Ext^i_{A \# G^*}(M \otimes_KKG^*,(A/J) \# G^*) \times \Ext^j_{A \#
G^*}((A/J) \# G^*,(A/J) \# G^*) & \stackrel{\nu}{\rightarrow} &
\Ext^{i+j}_{A \# G^*}(M \otimes_KKG^*,(A/J) \# G^*)
\end{array}$$} where $\nu$ is Yoneda product and $\mu$ is defined
by $\mu(\bar{\az} \otimes p_g, \bar{\bz} \otimes
p_h):=\bar{\az}(\bar{\bz})_{gh^{-1}} \otimes p_h$.

Let $\mathbb{Q}=(Q_k, \phi_k)_{k \ge 0}$ be a finitely generated
graded projective resolution of $A/J$. Let $\bar{\az} \in
\Ext_A^i(M, A/J)$ and $\bar{\bz} \in \Ext_A^j(A/J, A/J)$ with
$\alpha \in \Hom_A(\Omega^i(M), A/J)$ or $\Hom(P_i, A/J)$ and $\beta
\in \Hom_A(\Omega^j(A/J), A/J)$ or \linebreak $\Hom(Q_j, A/J)$. Then
$\bar{\az}(\bar{\bz})_{gh^{-1}}=\overline{\beta_{gh^{-1}}\Omega^j(\alpha)}$.
Since $\mathbb{P}=(P_k, \psi_k)_{k \ge 0}$ is a graded projective
resolution of $M$, $\mathbb{P}(i)=(P_{k+i}, \psi_{k+i})_{k\ge 0}$ is
a graded projective resolution $\Omega^i(M)$ for all $i\ge 0$.
Furthermore $\mathbb{P}(i) \otimes_KKG^*=(P_{k+i} \otimes_KKG^*,
\psi_{k+i} \otimes 1)_{k\ge 0}$ is an $A \# G^*$-module projective
resolution of $\Omega^i(M)\otimes_KKG^*$ for all $i\ge 0$. Suppose
$\alpha: \Omega^i(M)\rightarrow A/J$ induces a chain map
$(\alpha_j)_{j\ge 0}: \mathbb{P}(i)\rightarrow \mathbb{Q}$. Then
$(\theta_{P_{j+i}}(\az_j \otimes p_g))_{j \ge 0}: \mathbb{P}(i)
\otimes_KKG^* \rightarrow \mathbb{Q} \otimes_KKG^*$ is exactly a
chain map induced by the map $\theta_{P_i}(\alpha \otimes p_g)$.
Therefore $\theta_{P_{j+i}}(\Omega^j(\alpha)\otimes
p_g)=\Omega^j(\theta_{P_i}(\alpha \otimes p_g))$.

For any $x \in (P_{j+i})_l$ and $p_m \in KG^*$, we have
$$\begin{array}{lll} &&\theta_{Q_j}(\beta \otimes
p_h) \Omega^j(\theta_{P_i}(\az \otimes p_g))(x \otimes p_m)\\
& = & \theta_{Q_j}(\beta \otimes
p_h) \theta_{P_{j+i}}(\Omega^j(\alpha) \otimes p_g))(x \otimes p_m)\\
& = & \theta_{Q_j}(\beta \otimes
p_h)((\Omega^j(\alpha)(x))_{mlg^{-1}} \otimes p_g)\\
& = & \beta((\Omega^j(\alpha)(x))_{mlg^{-1}})_{mlh^{-1}}\otimes
p_h\\ & = & (\beta_{gh^{-1}}\Omega^j(\alpha)(x))_{mlh^{-1}} \otimes
p_h\\ & = & \theta_{P_{j+i}}(\beta_{gh^{-1}}\Omega^j(\alpha) \otimes
p_h)(x \otimes p_m), \end{array}$$

Thus $\bar{\theta}_{i+j}((\bar{\az} \otimes p_g)(\bar{\bz} \otimes
p_h))=\bar{\theta}_{i+j}(\bar{\az}(\bar{\bz})_{g^{-1}h} \otimes
p_h)=\bar{\theta}_{i+j}(\overline{\beta_{gh^{-1}}\Omega^j(\alpha)}
\otimes p_h)= \overline{\theta_{P_{j+i}}
(\beta_{gh^{-1}}\Omega^j(\alpha) \otimes
p_h)}=\overline{\theta_{Q_j}(\bz \otimes
p_h)\Omega^j(\theta_{P_i}(\az \otimes
p_g))}=\overline{\theta_{P_i}(\az \otimes p_g)}$ \linebreak
$\overline{\theta_{Q_j}(\bz \otimes p_h)}=\bar{\theta}_i(\bar{\az}
\otimes p_g)\bar{\theta}_j(\bar{\bz} \otimes p_h)$. So the diagram
above is commutative.

Let $J$ and $J_G$  be the Jacobson radical and graded Jacobson
radical of $A$ respectively. Since $\char K \nmid |G|$, by
\cite[Theorem 4.1 and Theorem 4.4]{CM}, we have $J(A \# G^*)=J_G \#
G^*=J \# G^*$. Moreover, $(A \# G^*)/J(A \# G^*) \cong (A/J)\# G^*$.

Therefore $M \otimes_KKG^*$ is a quasi-Koszul $A \# G^*$-module if
and only if $$\begin{array}{lll} &&\Ext^{i+j}_{A \# G^*}(M
\otimes_KKG^*,(A/J) \# G^*)\\ & = & \Ext^i_{A \# G^*}(M
\otimes_KKG^*,(A/J) \# G^*)\Ext^j_{A \# G^*}((A/J) \# G^*,(A/J) \#
G^*) \end{array}$$ for all $i,j \geq 0$, if and only if
$$\begin{array}{lll} && \Ext^{i+j}_A (M ,A/J) \otimes _KKG^*\\
& = & (\Ext^i_A (M ,A/J) \otimes _KKG^*)(\Ext^j_A (A/J,A/J) \otimes
_KKG^*)\\ & = & (\Ext^i_A (M,A/J) \Ext^j_A (A/J,A/J)) \otimes _KKG^*
\end{array}$$ for all $i,j \geq 0$, if and only if $$\Ext^{i+j}_A (M,A/J) =
\Ext^i_A (M,A/J) \Ext^j_A (A/J,A/J)$$ for all $i, j \geq 0$, if and
only if $M$ is a quasi-Koszul $A$-module. In particular, taking
$M=A/J$ we have the algebra $A \# G^*$ is quasi-Koszul if and only
if so is $A$. Moreover, according to the commutative diagram above,
the map $(\bar{\theta}_i)_{i \geq 0}$ is a positively graded
$K$-algebra isomorphism $\Ext^*_A(A/J,A/J) \# G^* \cong \Ext^*_{A \#
G^*}((A \# G^*)/J(A \# G^*), (A \# G^*)/J(A \# G^*))$.
\hfill{$\Box$}

\begin{remark} Let $A$ be a finite-dimensional
$G$-graded $K$-algebra, where $G$ is a finite group satisfying
$\char K \nmid |G|$. We can show in another way that $A$ is
quasi-Koszul if and only if so is $A \# G^*$: It follows from
Corollary~\ref{matrixcorollary}, Cohen-Montgomery's duality theorem
for coaction (cf. \cite[Theorem 3.5]{CM}) and
Theorem~\ref{skewtheorem} that $A$ is quasi-Koszul if and only if so
is $M_n(A) = (A \# G^*) \ast G$, if and only if so is $A \# G^*$.
\end{remark}

\section{Finite Galois covering and Koszulity} \label{galoiscovering}

In this section, we show that if a finite-dimensional connected
quiver algebra is Koszul then so are its Galois covering algebras
with finite Galois group $G$ satisfying $\char K \nmid |G|$. Indeed,
we have two ways to observe the relations between finite Galois
covering and (quasi-)Koszulity: One is via covering functor and
galois extension (cf. Section~\ref{galoisextension}), the other is
via covering of quiver with relations (cf.
Section~\ref{smashcovering}).

\subsection{Galois extension, covering functor and Koszulity} \label{galoisextension}

Let $K$ be an algebraically closed field. Let $A'$ be a $K$-algebra
and $G$ a finite group acting on $A'$ as a group of
$K$-automorphisms. We say that the pair $(A',G)$ is {\it pregalois}
if $A'$ is a projective $A' \ast G$-generator (cf. \cite[Section
1]{ARS1989}). We say that the pregalois pair $(A',G)$ is {\it left}
(resp. {\it right}) {\it galois} if $A'^G/\ann_{A'^G}S$ is a
semisimple artinian $K$-algebra for each simple left (resp. right)
$A'$-module $S$. We say that the pair $(A',G)$ is {\it galois} if it
is both left and right galois (cf. \cite[Section 4]{ARS1989}). In
this case, the $K$-algebra $A'$ is called a {\it galois extension}
of $A:=A'^G$. In case $A'$ is a finite-dimensional basic
$K$-algebra, $(A',G)$ is galois if and only if the induced action of
$G$ on the isomorphism classes of simple $A'$-modules is free (cf.
\cite[Proposition 5.3]{ARS1989}).

A {\it $K$-category} $\underline{A}$ is a preadditive category in
which the morphism sets are $K$-vector spaces and the compositions
are $K$-bilinear. A $K$-category $\underline{A}$ is called a {\it
locally bounded $K$-category} if it satisfies: (1) for each object
$a$ in $\underline{A}$, $\End_{\underline{A}}(a)$ is a local ring;
(2) for each pair of objects $a$ and $b$ in $\underline{A}$,
$\dim_K\Hom_{\underline{A}}(a,b)< \infty$; (3) distinct objects of
$\underline{A}$ are nonisomorphic; and (4) for each $a$ in
$\underline{A}$ there are only a finite number of objects $b$ in
$\underline{A}$ such that $\Hom_{\underline{A}}(a,b) \neq 0$ or
$\Hom_{\underline{A}}(b,a) \neq 0$ (cf. \cite[Section 2.1]{BG}).

If $\underline{A}'$ and $\underline{A}$ are locally bounded
$K$-categories then a $K$-linear functor $\underline{F}:
\underline{A}' \rightarrow \underline{A}$ is called a {\it
covering functor} if: (1) $\underline{F}$ is surjective on
objects; and (2) for each $a'$ in $\underline{A}'$ and $a$ in
$\underline{A}$, $\underline{F}$ induces isomorphisms $\coprod_{b'
\in \underline{F}^{-1}(a)}\Hom_{\underline{A}'}(b',a') \rightarrow
\Hom_{\underline{A}}(a,\underline{F}(a'))$ and $\coprod_{b' \in
\underline{F}^{-1}(a)}\Hom_{\underline{A}'}(a',b') \rightarrow
\Hom_{\underline{A}}(\underline{F}(a'),a)$ (cf. \cite[Section
3.1]{BG}).

Let $\underline{A}'$ be a finite locally bounded $K$-category.
Assume that $G$ is a finite group of $K$-automorphisms of
$\underline{A}'$ such that the induced action on the  objects is
free. Then, by \cite[3.1 Proposition]{Gabriel}, the quotient
category $\underline{A}'/G$ exists and the canonical projection
$\underline{A}' \rightarrow \underline{A}'/G$ is a covering functor,
which is called a {\it Galois covering with finite Galois group
$G$.}

Let $A$ be a finite-dimensional basic $K$-algebra. Then $A$ is
isomorphic to $KQ/I$ where $Q$ is a finite quiver and $I$ is an
admissible ideal of the path algebra $KQ$ (cf. \cite[Corollary
1.10]{ARS}). So $A$ can be viewed as a finite locally bounded
$K$-category $\underline{A}$ which is the quotient category of the
path category $KQ$ (cf. \cite[Section 2.1]{BG}).

Assume that $\underline{F}: \underline{A}' \rightarrow
\underline{A}$ is a Galois covering with finite Galois group $G$,
where the finite locally bounded $K$-category $\underline{A}'$
corresponds to a finite-dimensional basic $K$-algebra $A'=KQ'/I'$.
Then $(A',G)$ is a galois extension and $A=A'^G$ (cf. \cite[Theorem
6.2]{ARS1989}). It follows from Theorem~\ref{skewtheorem} that, in
the case of $\char K \nmid |G|$, $A'$ is quasi-Koszul if and only if
so is $A' \ast G$. By \cite[Propositon 1.2]{ARS1989}, we know the
category $\Mod A' \ast G$ is equivalent to the category $\Mod A'^G$,
i.e., $\Mod A$. Applying Theorem~\ref{moritaequivalencetheorem} we
have $A' \ast G$ is quasi-Koszul if and only if so is $A=A'^G$.
Thus, in the case of $\char K \nmid |G|$, $A'$ is quasi-Koszul if
and only if so is $A$.

If both $A=KQ/I$ and $A'=KQ'/I'$ are graded quiver algebras then $A$
is Koszul if and only if so is $A'$.

\begin{remark} Though we can obtain the relation between Koszulity
and finite Galois covering via galois extension, the underlying
field $K$ has to be assumed to be algebraically closed.
\end{remark}

\subsection{Covering of quiver with relations and Koszulity} \label{smashcovering}

Let $K$ be a field and $Q$ a quiver. From now on we always denote by
$Q_0$ (resp. $Q_1$) the vertex (resp. arrow) set of $Q$. For a path
$p$ in $Q$, we denote by $i(p)$ (resp. $t(p)$) the initial point
(resp. terminal) point of $p$. Let $Q'$ and $Q$ be quivers. Let $F:
Q' \rightarrow Q$ be a covering (in the topological sense) and $v'
\in Q'_0$. Let $\pi_1(Q',v')$ be the fundamental group of $Q'$ and
$F_*: \pi_1(Q',v') \rightarrow \pi_1(Q,F(v'))$ be the map induced by
$F$. We say a covering $F:Q' \rightarrow Q$ is {\it regular}  if
$F_*(\pi_1(Q',v'))$ is a normal subgroup of $\pi_1(Q,F(v'))$. In
this case, $Q$ is homeomorphic to the quotient space $Q'/G$ where $G
\cong \pi_1(Q,F(v'))/F_*(\pi_1(Q',v'))$ is the automorphism group of
the covering $F: Q' \rightarrow Q$ (cf. \cite[Chapter 5]{Massey}).

Some of the following contents are taken from \cite{G1983} for
reader's convenience. If $r = \sum^m_{i=1}k_ip_i$ is a $K$-linear
combination of paths in $Q$ and $v,w$ are vertices in $Q$, we denote
by $c_{v,w}(r)$ the {\it $(v,w)$-component of $r$} where $c_{v,w}(r)
: = \sum^n_{j=1}k_{i_j}p_{i_j}$ and $\{p_{i_j}\}^n_{j=1}$ is the
subset of $\{p_i\}^m_{i=1}$ of all those paths with initial point
$v$ and terminal point $w$.

A map $L: Q_0 \rightarrow Q'_0$ is called a {\it lifting} if $L(v)
\in F^{-1}(v)$ for all $v \in Q_0$. By uniqueness of path lifting,
if $p$ is a path with initial point $i(p)$ and terminal point $t(p)$
then we denote by $L(p)$ the unique path in $Q'$ with initial point
$L(i(p))$ such that $F(L(p))=p$. If $r = \sum^m_{i=1}k_ip_i$ is a
$K$-linear combination of paths in $Q$ we denote
$\sum^m_{i=1}k_iL(p_i)$ by $L(r)$.

We say $(Q, \rho)$ is a {\it quiver with relations} if $Q$ is a
quiver and $\rho$ is a set of $K$-linear combinations of paths in
$Q$.  Let $(Q',\rho')$ and $(Q,\rho)$ be quivers with relations. We
say $F: (Q',\rho') \rightarrow (Q,\rho)$ is a {\it morphism of
quivers with relations} if $F: Q' \rightarrow Q$ is a regular
covering of quivers satisfying: (1) $\rho' = \{L(r) | L : Q_0
\rightarrow Q'_0$ is a lifting and $r \in \rho \}$; (2) if $r' \in
\rho'$ and $v,w \in Q_0$ then there exist $v',w' \in Q'_0$ such that
$F(c_{v',w'}(r'))=c_{v,w}(F(r'))$.

Let $A=KQ/\langle \rho \rangle$ be a finite-dimensional connected
$K$-algebra where $\langle \rho \rangle$ denotes the admissible
ideal of $KQ$ generated by $\rho$. Let $F: (Q',\rho') \rightarrow
(Q,\rho)$ be a morphism of quivers with relations and finite group
$G$ its group of automorphisms. In this case, $F$ is called a {\it
Galois covering with finite Galois group} $G$. It follows from
\cite[Theorem 3.2]{G1983}) that there is a weight function $W: Q_1
\rightarrow G$ such that $A$ may be given the $G$-grading induced by
$W$. Moreover, $\mod KQ'/\langle \rho' \rangle$ is equivalent to the
category $\gr A$ of finite-dimensional $G$-graded $A$-modules. By
\cite[Theorem 2.2]{CM}, we have $\gr A$ is isomorphic to $\mod A \#
G^*$. Thus $\mod KQ'/\langle \rho' \rangle$ is equivalent to $\mod A
\# G^*$. It follows from Theorem~\ref{moritaequivalencetheorem} and
Theorem~\ref{smashtheorem} that, in the case of $\char K \nmid |G|$,
$KQ'/\langle \rho' \rangle$ is quasi-Koszul if and only if so is $A
\# G^*$, if and only if so is $A$.

Clearly, $\rho$ is a set of homogeneous relations in the length
grading if and only if so is $\rho'$. Thus, if $KQ/\langle \rho
\rangle$ is a graded quiver algebra with $\rho$ a set of homogeneous
relations in the length grading then so is $A'=KQ'/\langle \rho'
\rangle$. So we obtain the following theorem:

\begin{theorem} \label{koszulitycovering}
Let $KQ/\langle \rho \rangle$ be a finite-dimensional connected
Koszul $K$-algebra with $\rho$ a set of quadratic relations. Let $F:
(Q',\rho') \rightarrow (Q,\rho)$ be a Galois covering with finite
Galois group $G$ satisfying $\char K \nmid |G|$. Then $KQ'/\langle
\rho' \rangle$ is Koszul.
\end{theorem}

\begin{remark} Of course we also have the ``inverse'' of
Theorem~\ref{koszulitycovering} as follows: Let $KQ'/\langle \rho'
\rangle$ be a finite-dimensional connected Koszul $K$-algebra with
$\rho'$ a set of quadratic relations. Let $F: (Q',\rho') \rightarrow
(Q,\rho)$ be a Galois covering with finite Galois group $G$
satisfying $\char K \nmid |G|$. Then $KQ/\langle \rho \rangle$ is
Koszul. However, it is useless in practice.
\end{remark}

\section{Construction of Koszul algebras} \label{construction}

In this section, we provide a general construction of Koszul algebra
by Galois covering with finite cyclic Galois group. Moreover, as
examples, we construct Koszul algebras from exterior algebras and
Koszul preprojective algebras by finite Galois covering with either
cyclic or noncyclic Galois group.

\subsection{General construction} \label{generalconstruction}
Let $A=KQ/\langle \rho \rangle$ be a Koszul algebra with $\rho$ a
set of quadratic relations. Then one can construct a finite Galois
covering $F: (Q',\rho') \rightarrow (Q,\rho)$ with Galois group
$G=\mathbb{Z}_n, n \geq 2,$ satisfying $\char K \nmid n$, in the
following way: Let $Q'_0 := \{(v,\overline{r})|v \in
Q_0,\overline{r} \in \mathbb{Z}_n\}$, $Q'_1 := \{(a,\overline{r}):
(i(a),\overline{r}) \rightarrow (t(a),\overline{r+1}) |a \in Q_1, r
\in \mathbb{Z}\}$ and $\rho' = \{\sum^s_{i=1}k_i(a_{i2},
\overline{r+2})(a_{i1},\overline{r+1}) | \sum^s_{i=1}k_ia_{i2}a_{i1}
\in \rho, k_i \in K, a_{ij} \in Q_1, r \in \mathbb{Z}\}$. The
covering $F: (Q',\rho') \rightarrow (Q,\rho)$ is defined by
$(v,\overline{r}) \mapsto v$ and $(a,\overline{r}) \mapsto a$, which
is a Galois covering with Galois group $\mathbb{Z}_n$. By
Theorem~\ref{koszulitycovering}, the graded quiver algebras
$KQ'/\langle \rho' \rangle$ are Koszul for all $n \geq 2$ with
$\char K \nmid n$.

\subsection{Constructions from exterior algebras}

Exterior algebras is a quite important class of algebras, which
plays extremely important roles in many mathematical branches such
as algebraic geometry, commutative algebra, differential geometry.

Let $Q$ be the quiver given by one vertex $1$ and $m$-loops
$a_1,a_2,...,a_m$ with $m \geq 2$. Denote by $\rho$ the set
$\{a_i^2|1 \leq i \leq m\} \cup \{a_ia_j+a_ja_i| 1 \leq i<j \leq
m\}$. Then $A:=KQ/\langle \rho \rangle$ is the exterior algebra over
$K$ (cf. \cite{Mat}). It is well-known that $A$ is a Koszul algebra
and its quadratic dual is the algebra $K[x_1,...,x_m]$ of
polynomials in $m$ variables $x_1,...,x_m$.

Applying the general construction in
Section~\ref{generalconstruction}, we can obtain many new Koszul
algebras:

\begin{example} \label{cyclic 1} Define the quiver with relations
$(Q', \rho')$ by $Q'_0 := \{\overline{r}|\overline{r} \in
\mathbb{Z}_n\}$, $Q'_1 := \{(a_i,\overline{r}): \overline{r}
\rightarrow \overline{r+1} | 1 \leq i \leq m, r \in \mathbb{Z}\}$
and $\rho'=\{(a_i,\overline{r+1})(a_i,\overline{r})| 1 \leq i \leq
m, r \in \mathbb{Z}\} \cup
\{(a_i,\overline{r+1})(a_j,\overline{r})+(a_j,\overline{r+1})(a_i,\overline{r})|
1 \leq i<j \leq m, r \in \mathbb{Z}\}$. The covering $F : (Q',\rho')
\rightarrow (Q,\rho)$ is defined by $\overline{r} \mapsto 1$ and
$(a,\overline{r}) \mapsto a$, which is a Galois covering with Galois
group $\mathbb{Z}_n$. By Theorem~\ref{koszulitycovering}, the graded
quiver algebras $KQ'/\langle \rho' \rangle$ are Koszul for all $n
\geq 2$ with $\char K \nmid n$. \end{example}

We can construct more new Koszul algebras from exterior algebras by
finite Galois covering with cyclic Galois group:

\begin{example} \label{cyclic 2} For $1 \leq l \leq m$,
define the quiver with relations $(Q', \rho')$ by $Q'_0 :=
\{\overline{r}|\overline{r} \in \mathbb{Z}_n\}$, $Q'_1 :=
\{(a_i,\overline{r}): \overline{r} \rightarrow \overline{r+1} | 1
\leq i \leq l, r \in \mathbb{Z}\} \cup \{(a_i,\overline{r}):
\overline{r+1} \rightarrow \overline{r} | l+1  \leq i \leq m, r \in
\mathbb{Z}\}$ and $\rho'$ naturally induced by $\rho$, i.e.,
$\rho':=\{(a_i,\overline{r+1})(a_i,\overline{r}) | 1 \leq i \leq l,
r \in \mathbb{Z}\} \cup \{(a_i,\overline{r})(a_i,\overline{r+1}) |
l+1 \leq i \leq m,  r \in \mathbb{Z}\} \cup
\{(a_j,\overline{r+1})(a_i,\overline{r}) +
(a_i,\overline{r+1})(a_j,\overline{r}) | 1 \leq i < j \leq l, r \in
\mathbb{Z}\} \cup \{(a_j,\overline{r})(a_i,\overline{r+1}) +
(a_i,\overline{r})(a_j,\overline{r+1}) | l+1 \leq i < j \leq m, r
\in \mathbb{Z}\} \cup \{(a_j,\overline{r})(a_i,\overline{r}) +
(a_i,\overline{r-1})(a_j,\overline{r-1}) | 1 \leq i \leq l, l+1 \leq
j \leq m, r \in \mathbb{Z}\}$. The covering $F : (Q', \rho')
\rightarrow (Q,\rho)$ is defined by $\overline{r} \mapsto 1$ and
$(a_i,\overline{r}) \mapsto a_i$, which is a Galois covering with
Galois group $\mathbb{Z}_n$. By Theorem~\ref{koszulitycovering}, the
graded quiver algebras $KQ'/\langle \rho' \rangle$ are Koszul for
all $n \geq 2$ with $\char K \nmid n$.
\end{example}

\begin{remark} It is well-known that the exterior
algebras can be viewed as $\mathbb{Z}_2$-graded algebras. The graded
exterior algebras (called {\it Grassmann algebras} as well
\cite{CJK}) have applications in physics. In the case of $n=2$,
namely in the case the Galois group is $\mathbb{Z}_2$, by
\cite[Theorem 3.2]{G1983}, the category $\mod KQ'/\langle \rho'
\rangle$ is equivalent to the category $\gr KQ/\langle \rho \rangle$
of $\mathbb{Z}_2$-graded modules (or superrepresentations) over the
graded exterior algebra (or superalgebra) $KQ/\langle \rho \rangle$.
From this viewpoint, the algebras $KQ'/\langle \rho' \rangle$ in
Example~\ref{cyclic 1} and Example~\ref{cyclic 2} should be
important.
\end{remark}

We can also construct new Koszul algebras from exterior algebras by
finite Galois covering with non-cyclic Galois group:

\begin{example}  Consider the exterior algebra with $m=2$.
Define the quiver with relations $(Q',\rho')$ by $Q'_0 :=
\{(1,\overline{r})|\overline{r} \in \mathbb{Z}_n\} \cup \{(1',
\overline{r})|\overline{r} \in \mathbb{Z}_n\}$, $Q'_1 :=
\{(a_1,\overline{r}): (1,\overline{r}) \rightarrow (1',\overline{r})
| r \in \mathbb{Z}\} \cup \{(a'_1,\overline{r}): (1',\overline{r})
\rightarrow (1,\overline{r}) | r \in \mathbb{Z}\} \cup
\{(a_2,\overline{r}): (1,\overline{r}) \rightarrow
(1,\overline{r+1}) | r \in \mathbb{Z}\} \cup \{(a'_2,\overline{r}):
(1',\overline{r}) \rightarrow (1',\overline{r+1}) | r \in
\mathbb{Z}\}$ and $\rho'$ naturally induced by $\rho= \{a^2_1,
a^2_2, a_1a_2+a_2a_1\}$. The covering $F : (Q',\rho') \rightarrow
(Q,\rho)$ is defined by $(1,\overline{r}) \mapsto 1$,
$(1',\overline{r}) \mapsto 1$, $(a_i,\overline{r}) \mapsto a_i$ and
$(a'_i,\overline{r}) \mapsto a_i$, which is a Galois covering with
Galois group $\mathbb{Z}_2 \times \mathbb{Z}_n$. By
Theorem~\ref{koszulitycovering}, the graded quiver algebras
$KQ'/\langle \rho' \rangle$ are Koszul for all $n \geq 2$ with
$\char K \nmid 2n$.
\end{example}

We can also construct new Koszul algebras from exterior algebras by
finite Galois covering with Galois group not the direct product of
cyclic groups:

\begin{example}  Consider the exterior algebra with $m=2$.
Define the quiver with relations $(Q', \rho')$ by $Q'_0 :=
\{(1,\overline{r})|\overline{r} \in \mathbb{Z}_n\} \cup \{(1',
\overline{r})|\overline{r} \in \mathbb{Z}_n\}$, $Q'_1 :=
\{(a_1,\overline{r}): (1,\overline{r}) \rightarrow (1',\overline{r})
| r \in \mathbb{Z}\} \cup \{(a'_1,\overline{r}): (1',\overline{r})
\rightarrow (1,\overline{r}) | r \in \mathbb{Z}\} \cup
\{(a_2,\overline{r}): (1,\overline{r}) \rightarrow
(1,\overline{r+1}) | r \in \mathbb{Z}\} \cup \{(a'_2,\overline{r}):
(1',\overline{r+1}) \rightarrow (1',\overline{r}) | r \in
\mathbb{Z}\}$ and $\rho'$ naturally induced by $\rho= \{a^2_1,
a^2_2, a_1a_2+a_2a_1\}$. The covering $F : (Q',\rho') \rightarrow
(Q,\rho)$ is defined by $(1,\overline{r}) \mapsto 1$,
$(1',\overline{r}) \mapsto 1$, $(a_i,\overline{r}) \mapsto a_i$ and
$(a'_i,\overline{r}) \mapsto a_i$, which is a Galois covering with
Galois group the dihedral group $D_{2n}$. By
Theorem~\ref{koszulitycovering}, the graded quiver algebras
$KQ'/\langle \rho' \rangle$ are Koszul for all $n \geq 2$ with
$\char K \nmid 2n$.
\end{example}

\subsection{Constructions from Preprojective algebras}

Preprojective algebras is a very nice class of algebras, which
plays important roles in differential geometry and quantum groups.

Let $Q=(Q_0,Q_1)$ be a finite connected quiver without oriented
cycle. Then the double quiver $\overline{Q}$ is defined by
$\overline{Q}_0 := Q_0$ and $\overline{Q}_1:=Q_1 \cup \{a^*: t(a)
\rightarrow i(a) | a \in Q_1\}$. Suppose that $\overline{I}$ is the
ideal of $K\overline{Q}$ generated by all elements of the form
$\sum_{a \in Q_1}(aa^*-a^*a)$. Then the algebra ${\cal
P}(Q):=K\overline{Q}/\overline{I}$ is called the {\it preprojective
algebras of the quiver $Q$} (cf. \cite{R1998}).

Note that ${\cal P}(Q)$ is Koszul if and only if $Q$ is either
$\mathbb{A}_1, \mathbb{A}_2$ or is not a Dynkin diagram
\cite[Theorem 1.9]{Mp}. In the case of $Q=\mathbb{A}_1$, ${\cal
P}(\mathbb{A}_1)=K$. In the case of $Q=\mathbb{A}_2$, the algebra
${\cal P}(\mathbb{A}_2)$ is radical square zero, so is its finite
Galois covering algebras. Thus we cannot construct new Koszul
algebras from ${\cal P}(\mathbb{A}_1)$ and ${\cal P}(\mathbb{A}_2)$
using our approach. Since ${\cal P}(Q)$ is finite-dimensional if and
only if $Q$ is a Dynkin quiver, our approach cannot be applied to
preprojective algebras directly. However, we can construct many new
Koszul algebras from the quadratic dual of preprojective algebras,
which are finite-dimensional.

Let $K$ be an algebraically closed field. It follows from
\cite[Theorem 1.9]{Mp} that, in case $Q$ is a tree but not a Dynkin
diagram, the quadratic dual of ${\cal P}(Q)$ is just the trivial
extension ${\cal D}(Q)$ of the path algebra $K\tilde{Q}$, where
$\tilde{Q}$ is of the same underlying graph as $Q$ and each vertex
is either a source or a sink. By \cite[Theorem 1.8]{Mp}, ${\cal
D}(Q) \cong K\hat{Q}/\hat{I}$ where $\hat{Q}$ is defined by
$\hat{Q}_0:=Q_0$ and $\hat{Q}_1:=\overline{Q}_1$ and the ideal
$\hat{I}$ is generated by the set $\hat{\rho}$ which consists of the
following relations: (1) $a^*a-b^*b$, if $i(a)=i(b)$ is a source in
$\hat{Q}$; (2) $aa^*-bb^*$, if $t(a)=t(b)$ is a sink in $\hat{Q}$;
(3) $b^*a$, if $t(a)=t(b)$ and $a \neq b$; (4) $ab^*$, if
$i(a)=i(b)$ and $a \neq b$.

On one hand, applying the general construction in
Section~\ref{generalconstruction} to $K\hat{Q}/\langle \hat{\rho}
\rangle$, we can construct many new Koszul algebras. On the other
hand, we can also construct more new Koszul algebras in the
following way:

\begin{example} Define the quiver with relations $(Q', \rho')$ by
$Q'_0 := \{(v,\overline{r})|v \in Q_0,\overline{r} \in
\mathbb{Z}_n\}$, $Q'_1 := \{(a,\overline{r}): (i(a),\overline{r})
\rightarrow (t(a),\overline{r}) | a \in Q_1, \overline{r} \in
\mathbb{Z}_n\} \cup \{(a^*,\overline{r}): (i(a^*),\overline{r})
\rightarrow (t(a^*),\overline{r+1}) | a \in Q_1, r \in \mathbb{Z}\}$
and $\rho'$ naturally induced by $\hat{\rho}$. The covering $F:
(Q',\rho') \rightarrow (\hat{Q},\hat{\rho})$ defined by
$(v,\overline{r}) \mapsto v$, $(a,\overline{r}) \mapsto a$ and
$(a^*,\overline{r}) \mapsto a^*$, is a Galois covering with Galois
group $\mathbb{Z}_n$. By Theorem~\ref{koszulitycovering}, the graded
quiver algebras $KQ'/\langle \rho' \rangle$ are Koszul for all $n
\geq 2$ with $\char K \nmid n$.
\end{example}

\end{document}